\documentclass[graybox,table]{svmult}

\usepackage{type1cm}        % activate if the above 3 fonts are
\usepackage{makeidx}         % allows index generation
\usepackage{graphicx}        % standard LaTeX graphics tool
\usepackage{multicol}        % used for the two-column index
\usepackage[bottom]{footmisc}% places footnotes at page bottom

\usepackage{multirow}
\usepackage{xcolor}

\usepackage{pgfplots}
\pgfplotsset{compat=1.15}
\usetikzlibrary{arrows}
\usetikzlibrary[patterns]
\usetikzlibrary{shapes.geometric,decorations.markings}

\definecolor{ffqqqq}{rgb}{1,0,0}
\definecolor{wwwwww}{rgb}{0.4,0.4,0.4}

\usepackage{newtxtext}       % 
\usepackage[varvw]{newtxmath}       % selects Times Roman as basic font

\makeindex             % used for the subject index
\usepackage{algorithm}
\usepackage{algpseudocode}

\usepackage{soul}

\usepackage{todonotes}
\begin{document}

\setcounter{page}{1}

\title*{On the Unmapped Tent Pitching for the Heterogeneous Wave Equation}
\author{Marcella Bonazzoli\orcidID{0000-0002-0284-5643},\\ Gabriele Ciaramella\orcidID{0000-0002-5877-4426}, \\ Ilario Mazzieri\orcidID{0000-0003-4121-8092}}
%\authorrunning{Marcella Bonazzoli, Gabriele Ciaramella, and Ilario Mazzieri}
\institute{
Marcella Bonazzoli \at Inria, Unité de Mathématiques Appliqueés, ENSTA, Institut Polytechnique de Paris, 91120 Palaiseau, France, \email{marcella.bonazzoli@inria.fr}
\and  
Gabriele Ciaramella \at MOX Lab, Dipartimento di Matematica, Politecnico di Milano, Italy, \email{gabriele.ciaramella@polimi.it} 
\and 
Ilario Mazzieri \at MOX Lab, Dipartimento di Matematica, Politecnico di Milano, Italy, \email{ilario.mazzieri@polimi.it}
}
\maketitle

\abstract*{The Unmapped Tent Pitching (UTP) algorithm is a space–time domain decomposition method for the parallel solution of hyperbolic problems. It was originally introduced for the homogeneous one-dimensional wave equation in [Ciaramella, Gander, Mazzieri, 2024]. UTP is inspired by the Mapped Tent Pitching (MTP) algorithm  [Gopalakrishnan, Schöberl, Wintersteiger, 2017], which constructs the solution by iteratively building polytopal space–time subdomains, referred to as tents. In MTP, each physical tent is mapped onto a space–time rectangle, where local problems are solved before being mapped back to the original domain.
In contrast, UTP avoids the nonlinear and potentially singular mapping step by computing the solution directly on a physical space–time rectangle that contains the tent, at the expense of redundant computations in the region outside the tent. In this work, we investigate several strategies to extend UTP to heterogeneous media, where the wave propagation speed is piecewise constant over two subregions of the domain. Among the considered approaches, the most efficient in terms of computational time is the one employing space–time subdomains with identical spatial and temporal dimensions in both regions, determined by the maximum propagation speed. }

\section{Introduction}

The Unmapped Tent Pitching algorithm (UTP) is a novel space-time domain decomposition technique for {wave-propagation} problems, first introduced in \cite{procsUTP} for the parallel solution of the homogeneous one-dimensional wave equation and later extended in \cite{UTP3D} to three-dimensional space domains. It shares the idea of the Mapped Tent Pitching algorithm (MTP) \cite{MTP2017}, which computes the solution by iteratively constructing polytopal space-time subdomains, called \emph{tents}. 
Tents are pitched following the characteristic lines of the wave equation (CFL condition), so that the local problems can be solved exactly within each tent. 
More precisely, in MTP the {solution} is obtained by: (i) transforming the physical tents into Cartesian space-time cylinders, (ii) computing the local solutions in these transformed subdomains and (iii) mapping the solution back into the original tents. 

The idea of UTP is to avoid the {non-linear (and singular)} mapping process by computing the solution directly in physical cartesian cylinders containing the tents, but at the cost of redundant computations \cite{procsUTP}. 
The height of these space-time cylinders is determined by the CFL condition, as for the MTP algorithm, while the length is the size of the space subdomain. 
{However, it can be observed that when the space–time subdomains of the UTP and the mapped tents of the MTP coincide, the computational costs of the corresponding subproblems are of the same order, provided that the same space–time grid and time integrator are used. In the UTP, the chosen time integrator must satisfy a CFL condition, whereas in the MTP this restriction applies not only within the mapped tents but also through an additional limitation on the tent height imposed by the singularity of the mapping.
Furthermore, this singularity also leads to a reduction in the accuracy of the resulting numerical solution; see, e.g., \cite{Gopalakrishnan2020StructureAwareRK}.}

{In this paper, we study the UTP when applied to the heterogeneous problem}
\begin{equation}
\label{eq:hetwave}
  \begin{cases}
    \partial_{tt} u(x,t) = c(x)^2 \partial_{xx} u(x,t), & \text{for } (x,t) \in \Omega\times(0,T), \\
    u(x,0) = f(x), \text{ and } \partial_t u(x,0) = g(x), & \text{for } x \in \Omega,\\
    u(0,t) = u(L,t) = 0, & \text{for } t \in (0,T].
  \end{cases}
\end{equation} 
{Here, \(\Omega = (0,L)\), \(T>0\). We assume for simplicity that the wave speed is piecewise} constant with \(c(x) = c_1 > 0\) for \(x \in (0,L/2]\) and \(c(x) = c_2 > 0\) for \(x \in (L/2,L)\), and \(f, g\) sufficiently regular functions.  

{While in the homogeneous case the choice of the space-time UTP subdomains is rather natural, in this work we aim at answering the question:}
\begin{flushleft}
    \em What is the optimal space-time decomposition in the case of a heterogeneous wave speed?
\end{flushleft}
{Here, the optimality is understood in terms of computational time, and will be precisely defined in the following sections.}

{The article is structured as follows. In section \ref{sec:UTP:homo}, we recall the UTP for a homogeneous medium, and briefly discuss its properties. Section \ref{sec:UTP:hetero} is devoted to the study of UTP for the heterogeneous problem \eqref{eq:hetwave}. Conclusions are drawn in section \ref{sec:concl}.}

\section{UTP in the homogeneous case}\label{sec:UTP:homo}

Before extending the UTP to the heterogeneous case, let us recall the algorithm in the homogeneous setting, where \(c(x)=c=c_1=c_2\) and the slope of the characteristic lines is \(\pm 1/c\) (see also \cite[Algorithm 2 and Fig.~3]{procsUTP}). 
Consider a set \(\{x_j\}_{j=0}^{N-1} \subset \overline{\Omega}\) of \(N\) equispaced points \(x_j = j L/(N-1)\), \(j=0,\dots,N-1\), with \(N\) odd, and a decomposition of \(\Omega\) into \(N-2\) overlapping subintervals 
\begin{equation}
  \label{eq:Ij}
  I_j = (x_{j-1},x_{j+1}), \, j=1, \dots, N-2.
\end{equation}
By setting \(\mathcal{R} = \{1,3,5,\dots, N-2\}\) and \(\mathcal{B} = \{2,4,6,\dots, N-3\}\), we distinguish between the non-overlapping \emph{red} subintervals \(I_j, j \in \mathcal{R}\), and the non-overlapping \emph{black} subintervals \(I_j, j \in \mathcal{B}\). 
In the UTP space-time rectangular subdomains are pitched alternately on the red and black subintervals. Following the characteristic lines, we set \(H = \frac{L}{c(N-1)}\). {In the first iteration (\(k=1\)), the red rectangles \(\mathcal{T}_j^k\),  \(j \in \mathcal{R}\), of height \(H\) are pitched. In the next iterations (\(k>1\)), black and red rectangles \(\mathcal{T}_j^k\) of height \(2H\) are pitched alternately, with \(j \in \mathcal{B}\) when \(k\) is even and \(j \in \mathcal{R}\) when \(k\) is odd. 
%Let us call \(T_j^k\) the time interval in \(\mathcal{T}_j^k\) (see its expression in lines~\ref{line:k1} and \ref{line:kge1} of Algorithm~\ref{alg:UTPhom}).  
Thus, in each iteration \(k\), local problems are solved in parallel in the space-time rectangles 
\begin{equation*}
    \mathcal{T}_j^k := I_j \times (v_j^{k-2},v_j^{k}),
\end{equation*}
with $v_j^k = v_j^{k-2} +2 H$, \(j \in \mathcal{R}\) for \(k\) odd, and \(j \in \mathcal{B}\) for \(k\) even. 
Notice that for the first red iteration ($k=1$) the length in time is different: $\mathcal{T}_j^1 := I_j \times (0,H)$. The same can happen for the last iteration.}
{Thus, given an approximation $u^{k-1}$ in $\Omega \times (0,T)$, the new approximation $u^k$ is computed by first solving (in parallel) the subproblems 
\begin{equation}
  \label{eq:localPBhom}
  \begin{cases}
    \partial_{tt} u_j^k = c^2 \partial_{xx} u_j^k & \text{in } \mathcal{T}_j^k = I_j \times (v_j^{k-2},v_j^{k}), \\
    u_j^k = u^{k-1} & \text{on }  \{x_{j-1},x_{j+1}\} \times (v_j^{k-2},v_j^{k}),\\
    u_j^k = u^{k-1} & \text{on }  I_j \times \{ v_j^{k-2} \},\\
    \partial_t u_j^k = \partial_t u^{k-1} & \text{on }  I_j \times \{ v_j^{k-2} \},\\
  \end{cases}
\end{equation} 
%\begin{equation}
%  \label{eq:localPBhom}
%  \begin{cases}
%    \partial_{tt} u_j^k(x,t) = c^2 \partial_{xx} u_j^k(x,t) & \text{for } (x,t) \in \mathcal{T}_j^k, \\
%    u_1^k(0,t) = 0 & \text{for } t \in (v_1^{k-1},v_1^k),\\
%    u_{N-1}^k(L,t) = 0 & \text{for } t \in (v_{N-1}^{k-2},v_{N-1}^{k}),\\
%    u_j^k(x_{j-1},t) = u_{j-1}^{k-1}(x_{j-1},t) & \text{for } t \in (v_j^{k-2},v_j^{k}), j \ne 1,\\
%    u_j^k(x_{j+1},t) = u_{j+1}^{k-1}(x_{j+1},t) & \text{for } t \in T_j^k, j \ne N-2.
%  \end{cases}
%\end{equation} 
for $j \in \mathcal{R}$ for $k$ odd and $j \in \mathcal{B}$ for $k$ even,
and then setting $u^k = u_j^k$ in $\mathcal{T}_j^k$ for the solved $j$ subdomains and $u^k = u_j^{k-1}$ elsewhere.
Clearly, \textcolor{black}{we complement problem \eqref{eq:localPBhom} with Dirichlet boundary conditions, i.e. $u_1^k=0$ on $\{0\} \times (v_j^{k-2},v_j^{k})$, $u_{N-2}^k = 0$ on $\{L\} \times (v_j^{k-2},v_j^{k})$ for any $k$ odd. Moreover,} at the first two iterations ($k= 0$ and $k=1$), one considers the initial conditions \(u_j^1(x,0) = f(x)\) and \(\partial_t u_j^1(x,0) = g(x)\).
The overall UTP procedure is detailed in Algorithm~\ref{alg:UTPhom}.}
\begin{algorithm}[t]
  \caption{Unmapped Tent Pitching in the homogeneous case}\label{alg:UTPhom}
  \begin{algorithmic}[1]
      \Require The subintervals \(I_j \subset \Omega\) defined in \eqref{eq:Ij} and {an initial guess function $u^0$ in $\overline{\Omega \times (0,T)}$ such that $u^0 = f$ and $\partial_t u^0 = g$ on $\Omega \times \{0\}$.}
      \State {Set \(H = \frac{L}{c(N-1)}\), \(\mathcal{R} = \{1,3,5,\dots, N-2\}\), \(\mathcal{B} = \{2,4,6,\dots, N-3\}\).}
      \State {Set} \(k=1\) and \(v_j^0 = 0\) for all \(j=1,\dots,N\).
      \While{\(\exists \,j\in\{1,\dots,N-2\}\) such that \(v_j^{k-1} \ne T\)} 
        \State {Set $u^k = u^{k-1}$.}
        \State Set \(J_k = \mathcal{R}\) if \(k\) is odd, and \(J_k = \mathcal{B}\) if \(k\) is even. 
        \If{\(k=1\)}
          \State For all \(j\in J_k\), set \(v_j^k = \min(T,H)\) and pitch the rectangle \(\mathcal{T}_j^k = I_j \times (0,v_j^{k})\). \label{line:k1}
        \Else
          \State For all \(j\in J_k\), set \(v_j^k = \min(T,v_j^{k-2}+2H)\) and pitch \(\mathcal{T}_j^k = I_j \times (v_j^{k-2}, v_j^{k})\).\label{line:kge1}
        \EndIf
        \State {For \(j\in J_k\), solve problems~\eqref{eq:localPBhom} in parallel to get \(u_j^k\) in \(\mathcal{T}_j^k\).} %and extend \(u_j^k\) by \(u^0\) above \(\mathcal{T}_j^k\). 
        \State {Update $u^k = u_j^k$ in \(\mathcal{T}_j^k\) for all \(j\in J_k\).}
        \State Update \(k=k+1\).
      \EndWhile
  \end{algorithmic}
\end{algorithm}

%%%%%%%%%%%%%%%%%%%

\begin{figure}
\centering
\includegraphics[width=\textwidth]{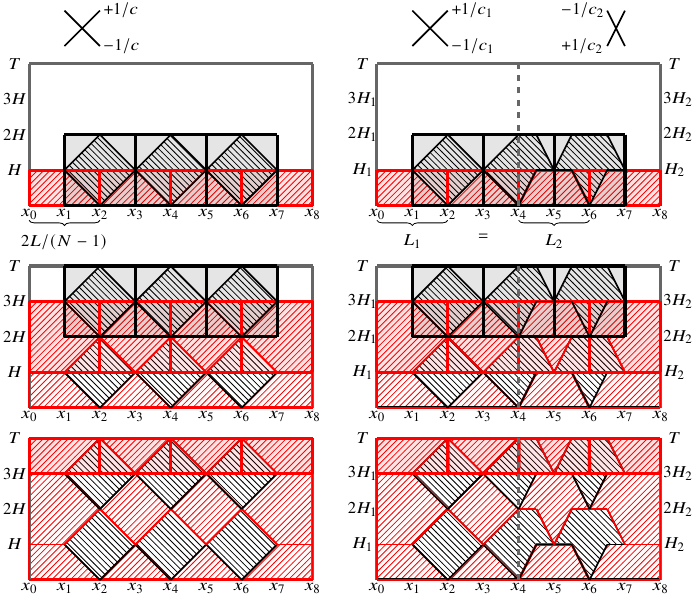}
\caption{Iterations of UTP with \(N=9\) in the homogeneous case (left), and in the heterogeneous case with the choice \(H_2 = H_1\), \(m_2=m_1=2\) (right). Here \(c_1=2c_2\). Red and black rectangles are the space-time subdomains constructed by UTP at odd and even iterations. We display iterations \(k=1, k=2\) at the top, iterations \(k=3, k=4\) in the middle, iteration \(k=5\) at the bottom of the figure. The hatched regions are the portions of the domain where UTP computes the exact solution.}
\label{fig:homANDcase1_H2eqH1_M2eqM1}
\end{figure}

{An example of UTP iterations is shown in Fig.~\ref{fig:homANDcase1_H2eqH1_M2eqM1} (left column).}
{The red and black subdomain solutions alternate. Notice that the red and black rectangles have height $2H = \frac{2L}{c(N-1)}$.}
As explained in \cite{procsUTP}, Algorithm 1 computes the exact solution below the tents (the hatched regions in Fig.~\ref{fig:homANDcase1_H2eqH1_M2eqM1}), and wrong approximations in the areas above the tents, which correspond to the regions where redundant computations are performed. However, note that the area of the rectangles (space-time subdomains) in the UTP scheme is the same as that of the MTP once the non-linear map is applied to the tent. For this reason, the local problem resolution {has the same computation cost for both algorithms.} The algorithm terminates when the exact solution is computed in the entire space-time domain \(\Omega\times(0,T)\).

%\section{Comparison of different UTP strategies in the heterogeneous case}\label{sec:UTP:hetero}
\section{{UTP in the heterogeneous case}}\label{sec:UTP:hetero}

We consider now the heterogeneous case \eqref{eq:hetwave} with \(c(x) = c_1 > 0\) for \(x \in (0,L/2]\) and \(c(x) = c_2 > 0\) for \(x \in (L/2,L)\). Without loss of generality, we assume that \(c_1 > c_2\), so the absolute value of the slope of the characteristic lines is smaller in the left region of the domain, i.e., \(|1/c_1| < |1/c_2|\). For simplicity, we assume \(T= L/(2c_1)\). 
{In contrast to the homogeneous case, the space-time subdomains in the two regions can have different lengths and heights. We denote by $L_j$ and $H_j$, for $j=1,2$, the subdomain lengths and heights in the two regions. 
Without loss of generality, we fix $L_1$ and $H_1$ and study the UTP behavior for varying $L_2$ and $H_2$.
In particular, we divide the left region \((0,L/2)\) into \(m_1\) non-overlapping subintervals (\(m_1 \in \mathbb{N}\backslash\{0\}\)) of equal length \(L_1 = L/(2m_1)\) and, following the characteristic lines, we set the height of the space-time rectangles (at the first iteration) to \(H_1 = L_1/(2c_1)\). In this way, the problem we aim to address in this work} is the following:
\begin{center} %flushleft
    \em Find the optimal values for the length \(L_2 = L/(2m_2)\), with \(m_2 \in \mathbb{N}\backslash\{0\}\), and the height \(H_2\) of the space-time rectangles for the right region \((L/2,L)\) such that the computational cost to {solve the problem in the entire $\Omega \times (0,T)$} is minimal.
\end{center} %flushleft
{Now, we assume that the computational cost of each subdomain problem is proportional to its space-time volume. This is the case, e.g., when an explicit time-stepping scheme is used. Furthermore, noting that at each iteration $m_1+m_2$ red subdomain problems or $m_1+m_2-1$ black subdomain problems are solved, we assume that $m_1+m_2$ parallel processes are available. Moreover, since at each iteration the subdomain problems are solved in parallel, we assume that the cost of  iteration is just biggest volume among those subdomains.}
%Note that with an explicit time scheme for the subdomain solves, the computational cost is proportional to the sum of the areas of the space-time rectangles used to cover \(\Omega\times(0,T)\), by considering, at each iteration \(k\), just the biggest among the red (or black) rectangles since the corresponding local problems are solved in parallel on \(m_1+m_2\) processes. 

\begin{figure}[t]
    \centering
    \includegraphics[width=\textwidth]{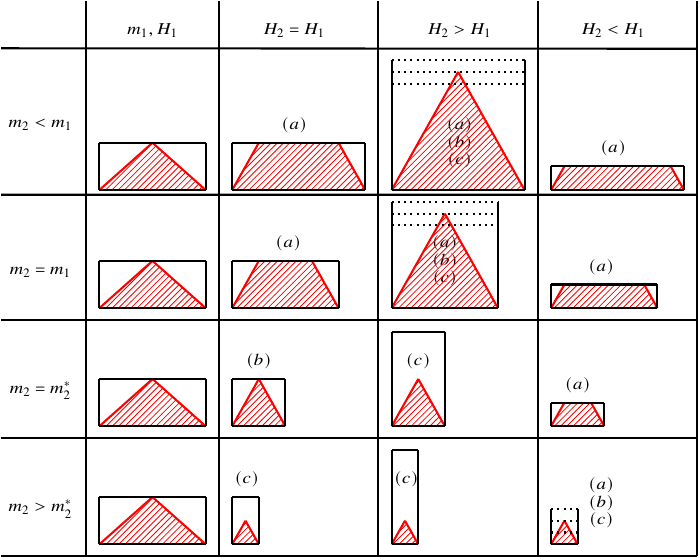}
    \caption{The possible choices for the length \(L_2 = L/(2m_2)\), with \(m_2 \in \mathbb{N}\backslash\{0\}\), and the height \(H_2\) of the space-time rectangles in the right region \((L/2,L)\), with respect to the reference rectangle of length \(L_1 = L/(2m_1)\) and height \(H_1 = L_1/(2c_1)\) in the left region \((0,L/2]\). We define \(m_2^*\) such that \(m_2^*/m_1 = c_1/c_2\), so \(m_2^* > m_1\) for the present case \(c_1>c_2\). The hatched region is the portion where the exact solution can be computed following the characteristic lines.}
    \label{fig:tabellina}
\end{figure}

{To find the optimal subdomains configuration, having fixed $L_1$ and $H_1$, it is necessary to identify all possible choices for \(m_2\) (hence $L_2$) and $H_2$. These are summarized in Fig.~\ref{fig:tabellina}.}
{Here, we show the space-time subdomains at the first iteration: the subdomains in the region with $c(x)=c_1$ in the first column of Fig.~\ref{fig:tabellina}, and those in the region with $c(x)=c_2$ in the second, third and fourth columns of Fig.~\ref{fig:tabellina}. Moreover, the different rows correspond to different sizes of the subdomain length $L_2$ given by $m_2$, while the different columns (second to third) correspond to different subdomain heights $H_2$.
For each case,} the hatched region is the portion in which the exact solution can be computed following the characteristic lines {if only exact initial data are available. 
The integer \(m_2^*\) is the one corresponding to \(L_2^* = L/(2 m_2^*)\) such that \(L_2^*/(2c_2) = L_1/(2c_1)\),} that is, the peaks of the `classical' triangular tents are reached at the same height in the two regions. We have \(m_2^*/m_1 = c_1/c_2\), so \(m_2^* > m_1\) for \(c_1>c_2\). 
{Moreover, the different choices are classified into three configurations: the peak of the triangular tent \((a)\) is not reached within the rectangular subdomain, \((b)\) is reached exactly at the top of the rectangular subdomain, and \((c)\) is reached strictly below the top of the subdomain.}
{Now, the goal is to compute the computational cost in all these cases, in order to identify the optimal configuration.}

{\underline{\(H_2 = H_1\) and \(m_2=m_1\)}. In this case,} the space-time rectangles have the same size in the left and right regions, with the height determined by \(c_1\), cf. Fig.~\ref{fig:tabellina}, second row and second column. 
As an example, with \(c_1=2c_2\) and \(m_2=m_1=2\), we show in Fig.~\ref{fig:homANDcase1_H2eqH1_M2eqM1} (right column) the \(5\) iterations needed for the UTP to compute the exact solution in the entire space-time domain \(\Omega\times(0,T)\), with \(T= L/(2c_1)\). More precisely, we display iterations \(k=1, k=2\) at the top, iterations \(k=3, k=4\) in the middle, iteration \(k=5\) at the bottom of the figure. With this first choice, UTP essentially proceeds as in the homogeneous case. The only difference is that in the heterogeneous case the portion where the exact solution is computed is not the same in the two regions of the domain, see the hatched portions in Fig.~\ref{fig:homANDcase1_H2eqH1_M2eqM1}. 
{Notice that Algorithm~\ref{alg:UTPhom} can be applied straightforwardly (with \(N = 2(2m_1)+1\)) by using \(c=c_1\) to pitch the rectangular subdomains and the true $c=c(x)$ for the local subdomain solves.} %{***Note that this is automatic when solving the local problems~\eqref{eq:localPBhom} by replacing \(c\) with \(c(x)\), and Algorithm~\ref{alg:UTPhom}, with \(c=c_1\) and \(N = 2(2m_1)+1\), needs no further modification.***CANCELLEREI QUESTA FRASE}

Let us now calculate the computational cost. %{***for the present strategy \(H_2 = H_1\), \(m_2=m_1\).CANCELLEREI}  
At each iteration \(k\), we use \(m_1+m_2 = 2m_1\) processes to solve in parallel the local problems on the rectangles: for \(k\) odd, we have \(2m_1\) red rectangles of the same size (the processes are perfectly balanced); for \(k\) even, we have \(2m_1 - 1\) black rectangles of the same size (one process is not used). The computational cost at \(k=1\) is proportional to the area \(A = L_1 H_1 = L_1^2/(2c_1) = L^2/(8 m_1^2 c_1)\). At the next iterations $k>1$ the cost is \(2A\), except the last one for which it is \(A\) again.  
To easily compute the total cost, we observe that the total height \(T= L/(2c_1)\) can be covered with \(T/H_1 = L/L_1 = 2m_1\) non overlapping red rectangles of area \(A\) (by considering the large ones as the union of two small ones), and with the same number \(2m_1\) of non overlapping black rectangles of area \(A\). Therefore, the total computational cost is 
\[
2(2m_1 A) = \frac{L^2}{2m_1c_1}.
\]
This is also shown in Fig.~\ref{fig:pipeline123} (b), where we report the pipeline of the computational work (time to solution) employed by the $m_1 + m_2$ parallel processes. It is easy to see that the computational workload is balanced and well-distributed among all the processes, with only one inactive process during the black iteration (the lighter-colored region at the bottom line).

%%%%%%%%%%%%%%%%%%%

\begin{figure}
\centering
\includegraphics[width=\textwidth]{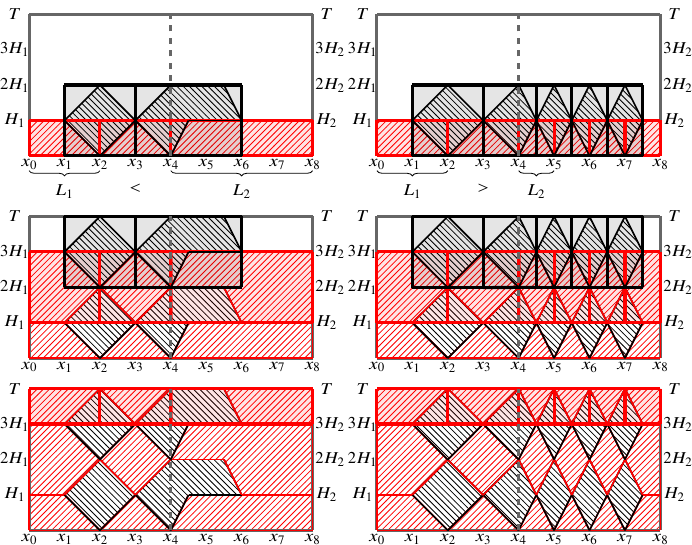}
\caption{Iterations of UTP with the choice \(H_2 = H_1\), \(1=m_2 < m_1=2\) (left), and the choice \(H_2 = H_1\), \(4=m_2=m_2^* \,(>m_1=2)\) (right).}
\label{fig:case2_H2eqH1_M2leM1ANDcase3_H2eqH1_M2eqM2star}
\end{figure}

%%%%%%%%%%%%%%%%%%%

{\underline{\(H_2 = H_1\) and \(m_2 < m_1\)}.}
We now consider the choice \(H_2 = H_1\), \(m_2 < m_1\), cf. Fig.~\ref{fig:tabellina} first row and second column, in which the rectangles are {larger} in the right region, see an example with \(m_1=2, m_2=1\) in Fig.~\ref{fig:case2_H2eqH1_M2leM1ANDcase3_H2eqH1_M2eqM2star} (left). By using again \(m_1+m_2\) processes, the local problems in the right part of the domain are more expensive to solve. {Hence,} \(A = L_2 H_1 = L_2 L_1/(2c_1) = L^2/(8 m_1 m_2 c_1)\), and the total computational cost is 
\[
2(2m_1 A) = \frac{L^2}{2m_2c_1} > \frac{L^2}{2m_1c_1}.
\]
Thus, this strategy is worse than the previous one {(\(H_2 = H_1\) and \(m_2 = m_1\)), as confirmed by comparing the corresponding pipelines in Fig.~\ref{fig:pipeline123} (a).}

{\underline{\(H_2 = H_1\) and \(m_2 = m_2^* > m_1\)}.}
An example {for this choice}, with \(m_1=2, m_2=m_2^*=4\), is displayed in Fig.~\ref{fig:case2_H2eqH1_M2leM1ANDcase3_H2eqH1_M2eqM2star} (right), cf. also Fig.~\ref{fig:tabellina} third row and second column. The {rectangular subdomains} are larger in the left region {and have area \(A = L_1 H_1\). Thus,} the total computational cost is the same as the first strategy ($H_2=H_1$ and $m_2=m_1$), but at the price of using more processes since \(m_1 + m_2> 2m_1\). {Moreover, no computational time improvement is achieved by increasing the number of processes, because \(m_2\) processes are partially inactive during the red or black iterations, cf. the pipeline in Fig.~\ref{fig:pipeline4} (a).}
%Moreover, using more processes yields no benefits in terms of computational time since half of them are inactive during either the red or black iterations, cf. Fig.~\ref{fig:pipeline123} (c).

%%%%%%%%%%%%%%%%%%%

\begin{figure}
\centering
\includegraphics[width=\textwidth]{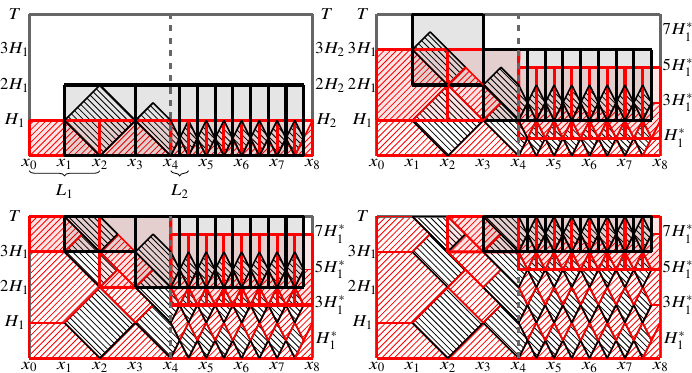}
\caption{First eight iterations of UTP with the choice \(H_2 = H_1\), \(8=m_2>m_2^*=4 \,(>m_1=2)\). We display iterations \(k=1, k=2\) at the top left, iterations \(k=3, k=4\) at the top right, iterations \(k=5, k=6\) at the bottom left, and iterations \(k=7, k=8\) at the bottom right of the figure.}
\label{fig:case4_H2eqH1_M2geM2star}
\end{figure}

%%%%%%%%%%%%%%%%%%%

%%%%%%%%%%%%%%%%%%%

\begin{figure}
\centering
\includegraphics[width=\textwidth]{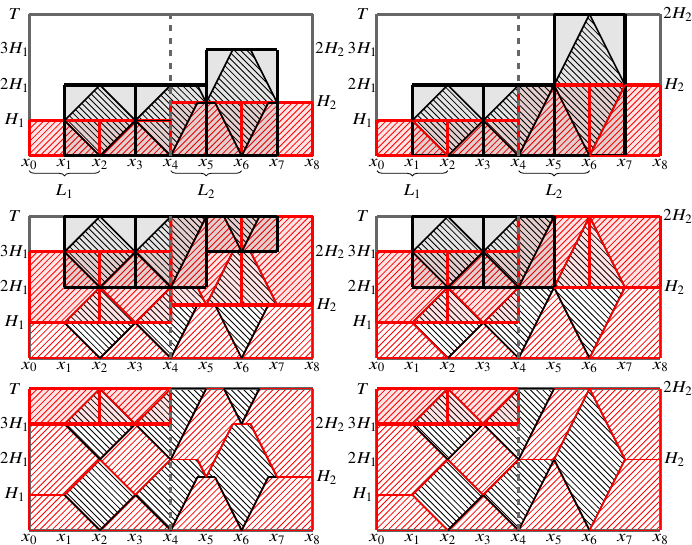}
\caption{Iterations of UTP with the choice \(H_2 > H_1\), \(m_2=m_1=2\), case \((a)\) (left), and the choice \(H_2 > H_1\), \(m_2=m_1=2\), case \((b)\) (right).}
\label{fig:case5_H2geH1_M2eqM1aANDcase5_H2geH1_M2eqM1b}
\end{figure}

%%%%%%%%%%%%%%%%%%%

{\underline{\(H_2 = H_1\) and \(m_2 > m_2^* > m_1\)}.}
{This case corresponds to Fig.~\ref{fig:tabellina} fourth row and second column, and it requires more iterations to compute the exact solution.} For instance, in Fig.~\ref{fig:case4_H2eqH1_M2geM2star}, we show the first eight iterations with \(m_1=2, m_2=8\), but nine iterations are actually needed, to be compared to five iterations in all the previous examples. Indeed, in the right region, since the portion where the exact solution is computed does not reach {the top of the rectangular subdomains}, the red rectangles at different iterations {need to} overlap (and similarly for the black ones). This leads to an exact solution advancing more slowly in the right region with multiples of an actual height $H_1^*<H_1$, which slows down the overall computation, as shown in the corresponding pipeline in Fig.~\ref{fig:pipeline4} (b). {Thus, the strategy \(H_2 = H_1\) and \(m_2 > m_2^* > m_1\) is outperformed by the previously discussed ones.}  

{\underline{\(H_2 > H_1\)}.}
This case corresponds to the third column in Fig.~\ref{fig:tabellina}, which could be regarded as a natural choice, since the characteristic lines in the right region are steeper. See two examples with \(m_2=m_1=2\) in Fig.~\ref{fig:case5_H2geH1_M2eqM1aANDcase5_H2geH1_M2eqM1b}. 
With \(H_2 > H_1\), the $m_2$ processes assigned to subdomains in the right region of the domain reach the final time \(T\) sooner than the $m_1$ ones in the left area. 
Nevertheless, this does not reduce the computational cost but rather causes communication delay within the processes, similarly to Fig.~\ref{fig:pipeline4} (b).

{\underline{\(H_2 < H_1\)}.} 
%Finally, choosing \(H_2 < H_1\) 
{This choice, corresponding to the fourth column in Fig.~\ref{fig:tabellina},  produces two effects.}
First, it conflicts with the characteristic slopes in the two regions. Moreover, it creates a symmetric issue to the one described above, that is, the final time $T$ is reached earlier in the left region than in the right one, resulting in a communication issue between processes. %\textcolor{red}{pipelines}.
%
%It not only goes against the slopes of the characteristic lines in the two regions, but also results in an issue which is symmetric to the one just above: the final time \(T\) is reached sooner in the left region than in the right region.
%
%

\bigskip
In conclusion:
\begin{flushleft}
    \em The optimal value for the length is $L_2 = L_1 = L/(2m_1)$, i.e., $m_2=m_1$, and for the height is \(H_2 = H_1 = L_1/(2c_1)\).
\end{flushleft}

{As discussed above, the resulting UTP algorithm for the heterogeneous wave equation~\eqref{eq:hetwave} is essentially  Algorithm~\ref{alg:UTPhom} using \(c=c_1\) (the maximum propagation speed) to build the rectangular subdomains, considering \(N = 2(2m_1)+1\), and solving, instead of \eqref{eq:localPBhom}, heterogeneous local problems: 
\begin{equation}
  \label{eq:localPBhet}
  \begin{cases}
    \partial_{tt} u_j^k(x,t) = c(x)^2 \partial_{xx} u_j^k(x,t) & \text{for $(x,t)$ in } \mathcal{T}_j^k = I_j \times (v_j^{k-2},v_j^{k}), \\
    u_j^k = u^{k-1} & \text{on }  \{x_{j-1},x_{j+1}\} \times (v_j^{k-2},v_j^{k}),\\
    u_j^k = u^{k-1} & \text{on }  I_j \times \{ v_j^{k-2} \},\\
    \partial_t u_j^k = \partial_t u^{k-1} & \text{on }  I_j \times \{ v_j^{k-2} \},\\
  \end{cases}
\end{equation} 
with the Dirichlet boundary conditions at \(x=0\) for \(u_1^k\), \(x=L\) for \(u_{N-2}^k\), for every \(k\) odd, and the initial conditions \(u_j^1(x,0) = f(x)\), \(\partial_t u_j^1(x,0) = g(x)\) for \(x \in I_j\) at the first iteration.}

% \begin{algorithm}[t]
%   \caption{Unmapped Tent Pitching in the heterogeneous case}\label{alg:UTPhet}
%   \begin{algorithmic}[1]
%       \Require Given \(2 m_1\) processes, the subintervals \(I_j \subset \Omega\) defined in \eqref{eq:Ij} with \(N = 4 m_1+1\) and an initial guess function \(u^0\). Set \(H = L/(4 m_1 c_1)\), \(\mathcal{R} = \{1,3,5,\dots, N-2\}\), \(\mathcal{B} = \{2,4,6,\dots, N-3\}\).
%       \State Let \(k=1\) and set \(v_j^0 = 0\) for all \(j=1,\dots,N\).
%       \While{\(\exists \,j\in\{1,...,N-2\}\) such that \(v_j^{k-1} \ne T\)} 
%         \State Set \(J_k = \mathcal{R}\) if \(k\) is odd, and \(J_k = \mathcal{B}\) if \(k\) is even. 
%         \If{\(k=1\)}
%           \State For all \(j\in J_k\), set \(v_j^k = \min(T,H)\) and pitch the rectangle \(\mathcal{T}_j^k = I_j \times (0,v_j^{k})\).
%         \Else
%           \State For all \(j\in J_k\), set \(v_j^k = \min(T,v_j^{k-2}+2H)\) and pitch \(\mathcal{T}_j^k = I_j \times (v_j^{k-2}, v_j^{k})\).
%         \EndIf
%         \State For \(j\in J_k\), solve problems~\eqref{eq:localPBhet} in parallel to get \(u_j^k\) in \(\mathcal{T}_j^k\), and extend \(u_j^k\) by \(u^0\) above \(\mathcal{T}_j^k\). 
%         \State Update \(k=k+1\).
%       \EndWhile
%   \end{algorithmic}
% \end{algorithm}

\begin{figure}
\centering
\includegraphics{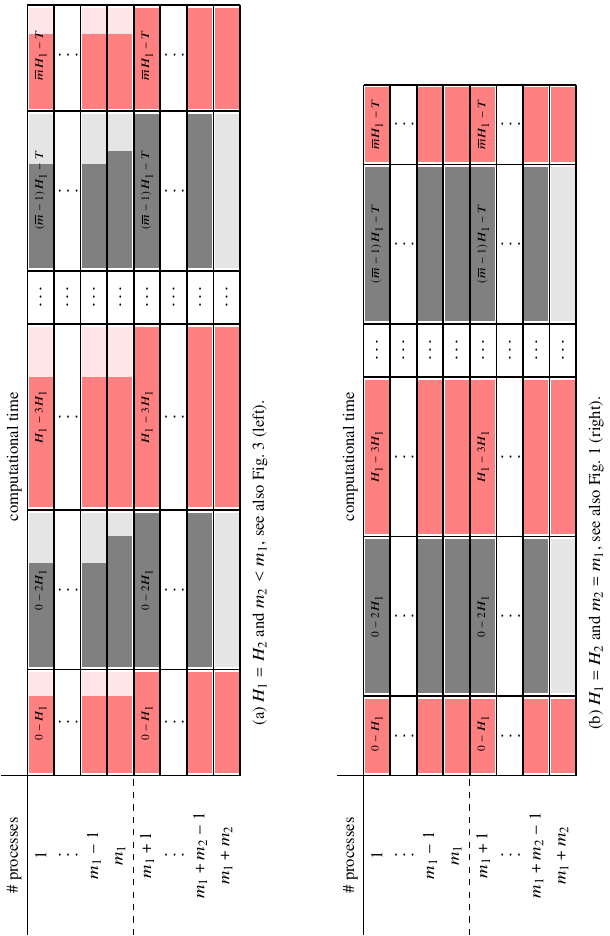}
\caption{Computational cost (time to solution) for different red-black domain decomposition strategies with $H_2 = H_1$, cf. Fig.~\ref{fig:tabellina} second column. The lighter-colored areas represent moments of inactivity of the corresponding process and $\overline{m} = 2m_1-1$.
}
\label{fig:pipeline123}
\end{figure}

\begin{figure}
\centering
\includegraphics{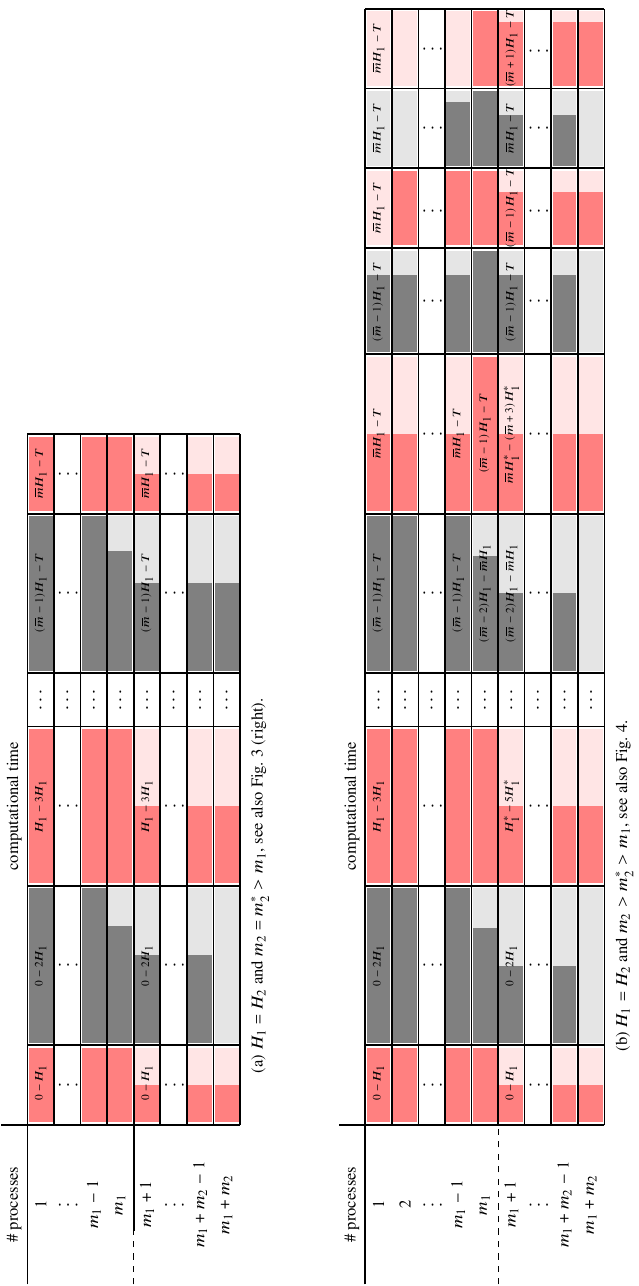}
\caption{Computational cost (time to solution) for different red-black domain decomposition strategies with $H_2 = H_1$, cf. Fig.~\ref{fig:tabellina} second column. The lighter-colored areas represent moments of inactivity of the corresponding process and $\overline{m}=2m_1-1$.}
\label{fig:pipeline4}
\end{figure}

\section{Conclusion}\label{sec:concl}
{In this work we studied the UTP for the wave equation (in one-dimensional space domains) with a piecewise constant propagation speed. In particular, the study focuses on analyzing the computational cost of different strategies for decomposing the space-time domain in order to identify the optimal one. We prove that the most efficient approach is the one employing space–time subdomains with identical spatial and temporal dimensions in the different material regions, determined by the maximum propagation speed. Future work will address the extension of these results to wave equations in two- and three-dimensional space domains.}

\begin{acknowledgement}
The work of G. Ciaramella and I. Mazzieri has been partially supported by the PRIN2022 grant ASTICE - CUP: D53D23005710006.
G. Ciaramella and I. Mazzieri are members of INdAM-GNCS group. The present research is part of the activities of ``Dipartimento di Eccellenza 2023-2027''.
%The authors would like to acknowledge the support of their respective institutions.
\end{acknowledgement}

\bibliographystyle{spmpsci}
\bibliography{references}

\end{document}